\documentclass{amsart}
\usepackage{amsmath}
\usepackage{amsthm}
\usepackage{amsfonts}
\usepackage{amssymb}
\usepackage[T1]{fontenc}
\usepackage{hyperref}
\usepackage[capitalise]{cleveref}
\usepackage{esint}
\usepackage{mathrsfs, mathtools, bm}
\usepackage{xcolor}
\usepackage{comment}
\usepackage{tikz}
\theoremstyle{plain}
\usetikzlibrary{calc,shapes.geometric,positioning,arrows.meta}
\usepackage{fullpage}
\usepackage{lmodern}
\usepackage{microtype}
\newtheorem{theorem}{Theorem}[section]
\newtheorem{defn}[theorem]{Definition}
\newtheorem{proposition}[theorem]{Proposition}

\newtheorem{lemma}[theorem]{Lemma}
\newtheorem{problem}[theorem]{Problem}

\newtheorem{remark}[theorem]{Remark}
\newcommand{\R}{\mathbb{R}}

\newcommand{\Sp}{\mathbb{S}}

\title{Uniqueness of $L_{p}$ Minkowski Problem in the supercritical range near the affine exponent}
\author{Devesh Rajpal}
\address{Mathematical Sciences Institute, Australian National University, Canberra ACT 2601, Australia}
\email{devesh.rajpal@anu.edu.au}
\date{}
\begin{document}

\begin{abstract}
    We study the isotropic \( L_{p} \)-Minkowski problem. The equation is affine invariant at the exponent \( p = -n-1 \), where all origin-centered ellipsoids with the volume of the unit ball are solutions, and it is known that the unit sphere is the only solution for \( p \in (-n-1, 1) \). We prove uniform \( C^{0} \) estimates for all solutions in the supercritical range \( p < -n-1 \). This allows us to prove a key spectral estimate of the Hilbert-Brunn-Minkowski operator on solutions. We then utilize the Rayleigh quotient method with position vectors as test functions for the operator adapted from Ivaki-Milman to prove uniqueness in the range \( p \in (-n-1- \delta_{n}, -n-1) \) for some \( \delta_{n} > 0 \). 
\end{abstract}
\maketitle
\section{Introduction}
    
    The \( L_{p} \)-Minkowski problem, introduced by Lutwak \cite{lutwak1993brunn} is a class of Monge-Amp\`{e}re problems on convex bodies that has attracted a lot of attention in the past two decades \cite{lutwak2004lp, hug2005lp, chou2006lp, BLYZ2013log, jian2015nonuniqueness, kolesnikov2022local}. It can be seen as an extension of the original Minkowski problem where we replace Minkowski sum by \( L_{p} \) Minkowski sum of the support functions. We describe the setup briefly.

    Let \( \Omega \subset \R^{n+1} \) be a convex body (convex domain with non-empty interior) with origin in the interior. We call the convex body smooth (or \( C^{k} \)) if the boundary hypersurface \( \mathcal{M} = \partial \Omega \) is smooth (or \( C^{k} \)) as an \( n \)-dimensional submanifold. Associated to any convex body is the support function \( \sigma_{\Omega} : \Sp^{n} \to \R \) which measures the distance from the origin to the supporting hyperplane in each direction.
    Define the operator \( A[f]_{ij} = \nabla_{i}\nabla_{j}f + g_{ij}f \), then the eigenvalues of \( A[\sigma] \) with respect to \( g_{ij} \) at \( z \) are the principal radii \( r_{1} = \kappa_{1}^{-1}, \dots, r_{n} = \kappa_{n}^{-1} \) of the hypersurface \( \partial \Omega \) at \( X(z) \). The Gaussian curvature is given by \( K(X(z))  = (\det_{g} A[\sigma])^{-1} = \det g_{ij} / \det A[\sigma]_{ij} \). 
    The Minkowski problem asks for the following: 
    \begin{problem}[Minkowski]
        Given a positive continuous function \( f \), does there exist a convex body \( \Omega \) such that \( K_{\Omega} = f^{-1} \)? If so, is it unique?
    \end{problem}
    One has to impose a necessary condition of Kazdan-Warner type on \( f \) coming from the divergence formula. The problem has been solved in affirmative in full regularity thanks to its resolution by combined work of Nirenberg \cite{nirenberg1953weyl}, Pogorelov \cite{pogorelov1978minkowski}, and Cheng-Yau \cite{cheng1976regularity}.  

    The support function formulation also allows an \( L_{p} \)-sum theory explored by Firey and Lutwak. 
    One can consider \( p \)-Minkowski sum and the corresponding Minkowski problem is the following. 
    \begin{problem}[\( L_{p} \)-Minkowski problem]
        Given a positive continuous function \( f \), does there exist a convex body \( \Omega \) with origin in the interior and support function \( \sigma \) such that \( \sigma^{1-p} \det A[\sigma] = f  \)? If so, is it unique?
    \end{problem}

    The \( L_{p} \)-Minkowski problem behaves differently for different values of \( p \). The supercritical range \( p \in (-\infty, - n-1) \) has no solutions for certain \( f \) \cite{chou2006lp,jian2015nonuniqueness,du2021nonuniqueness} and non-unique solutions for other \( f \), e.g. \( f \equiv 1 \). The equation is also connected to \( \alpha \)-Gauss curvature flow via the relation \( p = 1 - 1 / \alpha \). We will be concerned with the case where \( f \equiv 1 \), so solutions of the equation \begin{equation}
    \det_{g} A[\sigma] = \sigma^{p-1} \label{LpMinkowski}.
    \end{equation} 

    The above equation is Minkowski problem for \( p = 1 \), where uniqueness holds upto transformations. At \( p = n+1 \), the equation is scale invariant and spheres of all radii are solutions. For \( p > n+1 \), uniqueness follows from a direct application of maximum principle. Uniqueness for general \( f \), in the range \( p \in (1,n+1) \) is due to Lutwak \cite{lutwak1993brunn}.
    
    Assuming antipodal symmetry and regularity, Firey proved that the round sphere is the limiting shape of the Gauss curvature flow for surfaces \cite{firey1974shapes}. This corresponds to uniqueness of \cref{LpMinkowski} for \( p = 0 \). Andrews \cite{andrews1999gauss} proved it without the symmetry assumption for surfaces. For \( p \in (-n-1,1) \), that is \( \alpha > 1 / (n+2) \), Brendle-Choi-Daskalopoulos \cite{10.4310/ACTA.2017.v219.n1.a1} proved the uniqueness in all dimensions and the super-affine range \( p \in(-n-1,1) \). Ivaki-Milman \cite{ivaki2023uniqueness} later gave a short proof of uniqueness in the range \( p \in(-n-1,-1] \) based on a local Brunn-Minkowski inequality which is adapted in the present paper. At the affine exponent \( p = -n-1 \), the equation is invariant under \( \operatorname{SL}(n+1) \) transformations, and any origin-centred ellipsoid \( \mathcal{E} \) with volume constraint \( V(\mathcal{E}) = V(B^{n+1}) \) is a solution.

    In the supercritical range \( p < -n-1 \) uniqueness can fail. For curves \( (n = 1) \), Andrews \cite{andrews2003classification} classified all solutions: the circle is the only solution for \( - 7 \le p < -2 \), while for \( p < -7 \) there are also non-circular solutions with \( k \)-fold symmetry for each integer \( 3 \le k < \sqrt{2 - p} \). In higher dimensions, Andrews proved existence of non-spherical solutions \cite{andrews2000motion} with irreducible symmetry groups using Gauss curvature flow and a perturbation method. Du-Liu-Lu \cite{du2021nonuniqueness} also proved construction of non-spherical solutions with irreducible symmetry group using a calculus of variation argument for the associated functional. Andrews and the author's earlier work \cite{ABGCF} constructed an affine boosted Gauss curvature flow which was used to prove existence of non-spherical solutions with reducible symmetry groups (such as axial and dihedral symmetries). Hu-Ivaki \cite{hu2026centro} proved a centro-affine Poincar\'e inequality for the class of unconditional convex bodies (symmetric under reflection by coordinate planes) which proved the spectral estimate \cref{spectralestimate} for the unconditional class with \( \epsilon = 1 \) and uniqueness of \cref{LpMinkowski} for this symmetry class in the range \( p \in [-n-2,-n-1) \).

    Given the curves result, it is natural to ask whether a similar uniqueness results holds in higher dimensions for a range of \( p \) with \( p < -n-1 \). Our main result confirms this without any symmetry assumption.
    
    \begin{theorem} \label{maintheorem}
        There exists \( \delta_{n} >0  \) such that for every \( p \in (-n-1- \delta_{n}, - n -1) \) the unit sphere is the only solution of the isotropic \( L_{p} \)-Minkowski problem \cref{LpMinkowski}. 
    \end{theorem}
   
    The proof combines two results which may be of independent interest. The first is a \( C^{0} \) estimate that depends only on the dimension and is uniform on the whole supercritical range. 

    \begin{theorem}
        There exists \( C = C(n) > 1\) such that any solution of \cref{LpMinkowski} with \( p < -n-1 \) satisfies 
        \[ C^{-1} \le \sigma \le C \quad \text{ on }\, \Sp^{n} \quad \text{ and } \quad \frac{\rho_{+}(\Omega)}{\rho_{-}(\Omega)} \le C, \]
        where \( \rho_{+} \) and \( \rho_{-} \) denote the circumradius and inradius respectively. 
    \end{theorem}
    
    Du \cite[Theorem 1.5]{du2025uniqueness} proved \( C^{0} \) estimates that are uniform for \( p \) in compact subsets of \( (-\infty, - n -1) \). The novelty of the above theorem is that the constant doesn't degenerate as \( p \to -n-1 \), which allows us to prove uniform spectral estimate on solutions in \cref{spectralestimate}. The proof uses the results from the companion paper \cite{ABGCF} which proves convexity in the affine orbit of entropies \( Z_{1} \) and \( Z_{1-p} \), leading to entropy upper bounds. A lower entropy bound on solutions using Mahler volume ties up the complete result. 

    The second is a uniform spectral estimate of the Hilbert-Brunn-Minkowski (HBM) operator on solutions of \cref{LpMinkowski}. 

    \begin{theorem} \label{spectralestimate1}
        Let  \( q < -n-1 \) be fixed, \( \Omega \) be a solution of \cref{LpMinkowski} with \( q < p < -n-1 \) and \( \lambda_{2} \) be the second eigenvalue of HBM operator in \cref{HBM} associated to \( \Omega \). Then there exists \( \epsilon(n,q) > 0 \) such that \( \lambda_{2}(-\Delta_{\Omega}) \ge n+ \epsilon(n,q)  \).
    \end{theorem}

    While the current paper was in preparation, Wan \cite{wan2026spectralgapsuniquenessisotropic} proved the same result as \cref{maintheorem}, in an updated version of an earlier conditional uniqueness result. Liu-Lu \cite{liu2026remarksuniquenessnonuniquenessisotropic} also posted a preprint very recently proving it via a different method. \cref{maintheorem} was proven independently by the author and the method in the current paper builds on the earlier work \cite{ABGCF} of Andrews and the author to get the key spectral estimate, and the uniqueness proof is adapted from Ivaki-Milman \cite{ivaki2023uniqueness} and Du \cite{du2025uniqueness}.

    Being the last one to post their proof (as of yet), it gives the author a unique opportunity to do a quick comparison between different proofs of the same result about this very intersting problem. All three of the methods seem to crucially rely on isotropy lemma of \cref{isotropy} as it gives a lot of control on the geometric quantities associated to the solutions. The isotropy lemma implied isoperimetric estimates for isotropic bodies in \cite[Proposition 5.2]{ABGCF} which enabled convergence for the Affine Boosted Gauss Curvature Flow. Both Wan and Liu-Lu use polar body representation of the \( Z_{1-p} \) entropy of \cref{Z1-p} which is not surprising as the polar body representation is quite powerful and was crux of some of the geometric estimates in Andrews-Guan-Ni work on the convergence of \( \alpha \)-Gauss curvature flow \cite{andrews2016flow}. Wan's earlier work which was conditional on an isotropic position lemma can be related to author's work in \cite{ABGCF}. Liu-Lu's prove similar \( C^{0} \) estimates like in \cref{C0estimates} but circumvents Ivaki-Milman \cite{ivaki2023uniqueness} technique and go through a direct compactness argument without needing any spectral estimate. The current paper and Wan's result both rely on a spectral estimate to prove uniqueness by adapting Ivaki-Milman. Wan's impressive quantitative estimate of a more general spectral gap is done by building up on the work of Kolesnikov-Milman \cite{kolesnikov2022local} about the Neumann problem with Reilly formula and also relies on the recent resolution of the KLS conjecture.

    \subsection*{Acknowledgements}

    The author is grateful to his PhD supervisor Ben Andrews for many fruitful discussions on the subject. 

    \subsection*{AI disclosure.} AI assistance was used in part of lower bound estimate of entropy in \cref{entropyestimate}. The rest is author's own work. AI was also used for typo fixing and checking proofs.

    \section{Notation and Preliminary results}
    Let \( \Omega \subset \R^{n+1} \) be a convex body (convex domain with non-empty interior). We call the convex body smooth (or \( C^{k} \)) if the boundary hypersurface \( \mathcal{M} = \partial \Omega \) is smooth (or \( C^{k} \)) as an \( n \)-dimensional submanifold. Associated to any convex body is the support function \( \sigma_{\Omega} : \Sp^{n} \to \R \) defined by 
    \[ \sigma_{\Omega} (z) = \sup_{p \in \Omega} \left< z, p \right>. \]
    We will drop the subscript on \( \sigma \) when the convex body is clear from the context. A convex body \( \Omega \) is called locally uniformly convex if it is \( C^{2} \) and every point on the boundary has a unique normal. In this case, the Gauss map \( G : \partial \Omega \to \Sp^{n} \) is invertible with inverse given by 
    \[ X(z) =  \sigma(z)z + \nabla \sigma, \]
    where \( \nabla \) is the standard connection on \( \Sp^{n} \). The Gaussian curvature of the hypersurface at \( X(z) \) can be computed as a function of the support function. Let \( g_{ij} \) be the standard round metric on \( \Sp^{n} \) with respect to a local coordinate system. Define an operator \[ A[f]_{ij} = \nabla_{i}\nabla_{j}f + g_{ij}f, \] then the eigenvalues of \( A[\sigma] \) with respect to \( g_{ij} \) at \( z \) are the principal radii \( r_{1} = \kappa_{1}^{-1}, \dots, r_{n} = \kappa_{n}^{-1} \) of the hypersurface \( \partial \Omega \) at \( X(z) \). It follows that the Gaussian curvature is given by \( K(X(z))  = (\det_{g} A[\sigma])^{-1} = \det g_{ij} / \det A[\sigma]_{ij} \). We will often drop the subscript on \( \det A \) when it's clear from the context. The divergence formula can be used to express volume as a function of support function as 
    \[ V(\Omega) = \int_{\Omega}1 dx = \frac{1}{n+1}\int_{\Omega} \operatorname{div}(X)dx = \frac{1}{n+1}\int_{\partial \Omega}\left< X, \nu \right>d \mathcal{H}^{n} = \frac{1}{n+1}\int_{\Sp^{n}}\sigma \det A[\sigma]d \mu, \]
    where \( d \mu \) is the standard measure on \( \Sp^{n} \). The measure \( dV_{\Omega} = \frac{\sigma \det A[\sigma]}{n+1}d \mu \) associated to the volume functional is called the cone-volume measure of \( \Omega \).

    Let \( \Omega \) be a compact \( C^{2} \) convex body with strictly positive principal curvatures (so \( A[\sigma]_{ij} \) is invertible everywhere). The associated metric is given by \( g_{\Omega} = A[\sigma]_{ij}/ \sigma \). The Hilbert-Brunn-Minkowski (HBM) operator associated to \( \Omega \) is defined by 
    \begin{align}
        \Delta_{\Omega}f &= (A[\sigma]^{-1})^{ij} A[f \sigma]_{ij} - nf \\
        & = (A^{-1})^{ij}(\sigma\nabla^{2}_{ij}f + \nabla_{i}\sigma \nabla_{j}f + \nabla_{i}f \nabla_{j}\sigma). \label{HBM}
    \end{align}
    Milman \cite{milman2023centro} developed the theory for this operator and proved that it coincides with the centro-affine Laplacian of \( \Omega \), which explains the affine invariant properties of its spectrum. The HBM operator is self-adjoint with respect to the cone-volume measure \( dV_{\Omega} \). It can also be interpreted as the weighted Laplacian on the weighted manifold \( (\Sp^{n}, g_{\Omega}, dV_{\Omega}) \). The first eigenspace of \( \Delta_{\Omega} \) is spanned by linear functions adapted to \( \Omega \) with eigenvalue \( \lambda_{1}(\Delta_{\Omega}) = -n \); 
    \[ \Delta_{\Omega}\left( \frac{z^{i}}{\sigma} \right) = -n \frac{z^{i}}{\sigma}, \qquad i = 1,\dots,n+1. \]
    The Brunn-Minkowski inequality in its spectral formation, due to Hilbert is equivalent to the above fact. This also leads to a Poincar\'e inequality for the HBM operator.
    \begin{lemma}[Poincar\'e inequality]
        \cite[Lemma 8]{andrews1997monotone}
        Let \( f \in C^{2}(\Sp^{n}) \) with \( \int_{\Sp^{n}}f dV_{\Omega} = 0 \). Then \begin{equation}
            n\int_{\Sp^{n}}f^{2}dV_{\Omega} \le \int_{\Sp^{n}}f(-\Delta_{\Omega}f)dV_{\Omega} = \int_{\Sp^{n}}|\nabla f|_{g_{\Omega}}^{2}dV_{\Omega} 
        \end{equation}
        with equality if and only if \( f(z) = \frac{\left< z, v \right>}{\sigma(z)} \) for some \( v \in \R^{n+1} \). \label{Poincare}
    \end{lemma}
    
    Following the notation from \cite{ABGCF}, we work primarily with two scale-invariant functionals \( Z_{1} \) and \( Z_{1-p} \) defined for any \( C^{2} \) convex body as 
    \begin{equation}
        Z_{1}(\sigma) = \left( \frac{V(\sigma)}{V(B_{1})} \right)^{n(\frac{1}{1-p} -1)} \left( \fint_{\Sp^{n}}K^{ \frac{p}{1-p}}d \mu \right)^{n+1} = \left( \fint_{\Sp^{n}}\sigma \det A[\sigma]\right)^{ \frac{np}{1-p}} \left( \fint_{\Sp^{n}}\det A[\sigma]^{ \frac{-p}{1-p}}d \mu \right)^{n+1} 
    \end{equation}
    \begin{equation}
        Z_{1-p}(\sigma) =  \left( \frac{V(\sigma)}{V(B_{1})} \right)^{-p} \left( \fint_{\Sp^{n}} \sigma^{p} \right)^{n+1} = \left( \fint_{\Sp^{n}} \sigma \det A[\sigma] d \mu \right)^{-p} \left( \fint_{\Sp^{n}}\sigma^{p} \right)^{n+1} \label{Z1-p}
    \end{equation}
    which can be seen as two extremes; one expressed purely as a power of curvature term \( \det A \) and other purely as a power of the support function \( \sigma \). There is even a wider class of Orlicz-type functionals with solutions of \( L_{p} \)-Minkowski problem as critical points \cite{andrews1997monotone}.

    In \cite{ABGCF}, it was proved that the entropies \( Z_{\beta} \) for \( \beta > 0 \), which include \( Z_{1} \) and \( Z_{1-p} \) are convex in the affine orbit for \( p < -n-1 \). The solutions of \cref{LpMinkowski} are critical points of the \( Z_{1} \) functional, in particular they are critical points in their respective affine orbits. Spherical measures which are critical with respect to the \( \operatorname{SL}(n+1) \) action are called isotropic measures and are of particular interest for minimizing problems associated to geometric functionals such as area and quermassintegrals of convex bodies \cite{petty1961surface,giannopoulos1999isotropic,giannopoulos2000extremal,Lutwak_Yang_Zhang_2005}. 
    \begin{defn}
    A non-negative finite Borel measure \( \mu \) on \( \Sp^{n} \) is isotropic if 
    \[  \int_{\Sp^{n}}|x \cdot z|^{2}d \mu(z) = \frac{|\mu|}{n+1} \]
    for all \( x \in \Sp^{n} \) where \( |\mu| = \int_{\Sp^{n}}d \mu \). Equivalently, if for any linear map \( T :\R^{n+1} \to \R^{n+1}\), 
    \[ \int_{\Sp^{n}}\left< Tz,z \right>d \mu(z) =  \frac{\operatorname{tr}(T)|\mu|}{n+1}. \]
    \end{defn}

    By \cite[Proposition 3.3]{ABGCF}, the measure \( dZ_{1} = \det A[\sigma]^{ \frac{-p}{1-p}} d \mu \) is isotropic for solutions of \cref{LpMinkowski}. Note that this is same as isotropy of the cone-volume measure \( dV_{\Omega} = \frac{\sigma \det A[\sigma]}{n+1}d \mu \) for solutions. 
    \begin{proposition}
        Let \( \Omega \) be a solution of \cref{LpMinkowski}, then the cone-volume measure \( dV_{\Omega} \) is isotropic. In particular, 
        \[ \int_{\Sp^{n}}(z^{i})^{2}dV_{\Omega} = \frac{1}{n+1}\int_{\Sp^{n}}d V_{\Omega}. \] \label{isotropy}
    \end{proposition}
    Hu-Ivaki also proved isotropy of the cone-volume measure via a different method \cite{hu2025stability}.

    \section{Dimensional \(C^{0}\) estimates on a solution} \label{C0estimates}

    We turn to proving \( C^{0} \) estimates for a solution which are also closed near the endpoint \( p = -n-1 \). Earlier methods with \( Z_{1-p} \) entropy seem inadequate for this result \cite[Theorem 1.5]{du2025uniqueness}. The essential idea is to use isotropy as well as affine convexity in the supercritical range. As in the critical and subcritical range, it turns out that there are two-sided estimates for entropy/\( p \)-Blaschke-Santal\'o functional for a solution. 
   
    \begin{theorem} \label{entropyestimate}
        For a solution with \( p < -n-1 \) the entropy is bounded above and below by dimensional constants 
        \[ C_{1}(n) \le Z_{1}(\Omega) \le C_{2}(n). \]
        More explicitly, we have the following crude estimate for a solution of the \( L_{p} \)-Minkowski problem - 
        \[  (2(n+1))^{-(n+1)} \le \left( (2(n+1))^{-(n+1)} \right)^{\frac{-p}{1-p}} \le Z_{1}(\Omega) \le \left( (n+1)^{n+1} \right)^{ \frac{-p}{1-p}} \le (n+1)^{(n+1)}. \] 
    \end{theorem}
    \begin{proof}
        The upper bound follows from using affine minimizing property of the solutions with a John ellipsoid argument following \cite[Proposition 5.1]{ABGCF}. For the lower bound recall that the functionals \( Z_{1} \) and \( Z_{1-p} \) are related on a solution by the equation 
    \begin{equation}
        Z_{1}(\sigma) = Z_{1-p}(\sigma)^{\frac{1}{1-p}}. \label{Z1nZp}
    \end{equation}
    Recall that \( p = -n-1 \) is the affine exponent and the functional \( Z_{n+2} \) at this exponent is the Mahler volume product. Moreover, the \( p \)-Blaschke-Santal\'o functional for \( p < -n-1 \) is bounded by the Mahler product by the H\"older inequality as \begin{align}
        Z_{1-p}^{-\frac{1}{p}} & = \left( \fint_{\Sp^{n}}\sigma \det A[\sigma] d \mu \right) \left( \fint_{\Sp^{n}} \sigma^{p} d \mu\right)^{ \frac{n+1}{-p}}  \nonumber\\
        & \ge \left( \fint_{\Sp^{n}} \sigma \det A[\sigma]d \mu \right) \left( \fint_{\Sp^{n}}\sigma^{-n-1} d \mu \right) = V(\Omega)V(\Omega^{\circ})V(B_{1})^{-2}. \label{Holder}
    \end{align}
    
    Recall that the Mahler product is \( \operatorname{GL}(n+1) \) invariant. Let \( J \) be the John ellipsoid of \( \Omega \), so there exists a point \( q \in  \operatorname{Int}(\Omega)  \) such that \( J \subset \Omega - q \subset (n+1)J \). Let \( T \) be the positive-definite affine transformation such that \( TJ = B_{r}(0) \) is a ball of radius \( r \) about the origin. Then we have 
    \[ B_{r}(0) \subset T\Omega - Tq \subset B_{(n+1)r}.\]
    Since origin is contained in \( \operatorname{Int}(\Omega) \), it follows that \( |Tq| \le (n+1)r \), so \( T \Omega \subset B_{(n+1)r} + Tq \subset B_{2(n+1)r}  \). Taking polar duals, we get \( (B_{2(n+1)r}(0))^{\circ} = B_{\frac{1}{2(n+1)r}}(0) \subset (T \Omega)^{\circ} \). Combining these, we get a lower bound on the Mahler volume product, as 
    \begin{equation}
        V(\Omega)V(\Omega^{\circ}) = V(T \Omega)V((T \Omega)^{\circ}) \ge V(B_{r}(0))V(B_{\frac{1}{2(n+1)r}}(0)) = (2(n+1))^{-(n+1)}V(B_{1})^{2}. \label{Mahlerlb}
    \end{equation}
    Combining \cref{Z1nZp}, \cref{Holder}, and \cref{Mahlerlb} gives the desired result.
    \end{proof}
    \begin{remark}
        Mahler conjectured that the sharp lower bound for \( V(\Omega)V(\Omega^{\circ}) \) is \( \frac{(n+2)^{n+2}}{((n+1)!)^{2}} \). The crude bound \cref{Mahlerlb} suffices for our purposes.
    \end{remark}
    \begin{theorem}
        There exists a constant \( C(n) > 1 \) such that any solution of \cref{LpMinkowski} with \( p < -n-1 \) satisfies 
        \[ \frac{\rho_{+}(\sigma)}{\rho_{-}(\sigma)} \le C(n), \]
        where \( \rho_{+} \) and \( \rho_{-} \) denote the circumradius and inradius respectively.
    \end{theorem}
    \begin{proof}
    The solutions being \( Z_{1} \)-isotropic, satisfy the isoperimetric ratio from \cite[Proposition 5.2]{ABGCF}. Tracking the dependence on the parameter \( p \) leads to the following estimate: 
    \begin{equation}
        \frac{\rho_{+}(\sigma)}{\rho_{-}(\sigma)} \le C'(n)(n+1)^{\frac{(n+1)(p-1)}{p}}Z_{1}(\sigma)^{ \frac{1-p}{p}} \le C(n)Z_{1}(\sigma)^{ \frac{1-p}{p}}.\label{iso}
    \end{equation}
    The lower bound on \( Z_{1} \) establishes the isoperimetric estimate.
    \end{proof}
    
    \begin{proposition}
        There exists a constant \( C(n) > 0 \) such that any solution of \cref{LpMinkowski} with \( p \in (-\infty, -n-1) \) satisfies 
        \[ C(n)^{-1} \le \sigma(x) \le C(n), \qquad \forall x \in \Sp^{n}. \]
    \end{proposition}
    \begin{proof}
        For a solution, the normalized volume is given by 
        \begin{align*}
            \frac{V(\Omega)}{V(B^{n+1})} & = \fint_{\Sp^{n}}\sigma \det A[\sigma]d \mu = \fint_{\Sp^{n}}\det A[\sigma]^{ \frac{-p}{1-p}}d \mu.
        \end{align*}
        It follows that \begin{align*}
            \frac{V(\Omega)}{V(B^{n+1})} = (Z_{1}(\Omega))^{ \frac{1-p}{n+1-p}}
        \end{align*}
        is bounded above and below by \cref{entropyestimate}.
        Now by the isoperimetric estimate \( \rho_{-}(\Omega) \sim \rho_{+}(\Omega) \sim V(\Omega)^{\frac{1}{n+1}} \). This implies a dimensional upper bound on \( \rho_{+}(\Omega) \) and a dimensional lower bound on \( \rho_{-}(\Omega) \).
        Since origin is contained in the interior, we have \( \sigma(x) \le 2 \rho_{+} \le C(n) \). 

        For lower bound we adapt the method from \cite[Lemma 3.3]{du2025limiting}. The idea is to leverage convexity and upper bound of the support function to get a lower bound. Let \( M = \max \sigma \) and suppose that \( m = \min \sigma \) is attained at \( z_{0} \in \Sp^{n} \). As \( |X|^{2} = \sigma^{2} + |\nabla \sigma|^{2}\), using Lipschitz bound \( |\nabla \sigma| \le M \), we have
        \[ \sigma(z) \le m + M\theta, \qquad \forall z \in \Sp^{n}, \]
        where \( \theta = \operatorname{dist}(z,z_{0}) \) is the geodesic distance. Now, \begin{align*}
            C(n) \ge \fint_{\Sp^{n}}\sigma^{p}d \mu & \ge C_{1}(n) \int_{0}^{\pi /2} (m+M \theta)^{p}\theta^{n-1}d \theta \\
            & = C_{1}(n)M^{-n}m^{p+n} \int_{0}^{M \pi /2m} (1+t)^{p}t^{n-1}dt \\
            & \ge C_{1}(n)M^{-n}m^{p+n} \int_{0}^{1 /2} (1+t)^{p}t^{n-1}dt \\
            & \ge C_{2}(n)M^{-n}m^{p+n}2^{p}.
        \end{align*}
        Utilizing the dimensional upper bound we get a dimensional lower bound,
        \[ m \ge  \left(\frac{C_{2}(n)}{C(n)}\right)^{\frac{-1}{n+p}}M^{ \frac{n}{n+p}}2^{ \frac{-p}{p+n}}\ge C_{3}(n)2^{ \frac{-p}{p+n}} \ge C_{3}(n)\inf_{p \in (-\infty, -n-1)}2^{\frac{-p}{p+n}} \ge C_{4}(n), \]
        where we used \( p+n < -1 \) for the first inequality.
    \end{proof}
    \section{Spectral estimate of Hilbert-Brunn-Minkowski operator on a solution}

    In this section, we establish a spectral estimate of the HBM operator for a solution without any symmetry assumptions. The spectra of HBM operator has been a subject of research in the recent years because of its connection to log-Brunn-Minkowski conjecture \cite{kolesnikov2022local,milman2024sharp}. Its connection to \( \alpha \)-Gauss curvature flow is also apparent since it is the linearized operator of the flow. Let \( \Omega \) be a \( C^{2} \) locally uniformly convex closed convex body with origin in its interior. Following \cite{milman2023centro}, define the centro-affine Laplacian/HBM operator as 
    \begin{align*}
        \Delta_{\Omega}f &= (A[\sigma]^{-1})^{ij} A[f \sigma]_{ij} - nf \\
        & = (A^{-1})^{ij}(\sigma\nabla^{2}_{ij}f + \nabla_{i}\sigma \nabla_{j}f + \nabla_{i}f \nabla_{j}\sigma).
    \end{align*}
    \begin{theorem} \label{spectralestimate}
        Let  \( q < -n-1 \) be fixed and \( \Omega \) be a solution of the isotropic \( L_{p} \)-Minkowski problem \cref{LpMinkowski} with \( q < p < -n-1 \). Then there is an \( \epsilon(n,q) > 0 \) such that \( \lambda_{2}(-\Delta_{\Omega}) \ge n+ \epsilon(n,q)  \).
    \end{theorem}
    \begin{proof}
        By \( C^{0} \) estimates and the \( L_{p} \)-Minkowski equation we know that \( C(n,q)^{-1} \le \det A[\sigma] \le C(n,q) \). Now, \( X = \sigma z + \nabla \sigma \) and \( z \perp \nabla \sigma \), so \( |X|^{2} = \sigma^{2} + |\nabla \sigma|^{2} \). It follows that \( |\nabla \sigma| \le \sup |X| = \sigma_{\max} \)  and by \( C^{0} \) estimates, \( \sigma^{p-1} \) is uniformly H\"older. By Caffarelli's result \cite{caffarelli1989interior}, this implies that \( ||\sigma||_{C^{2, \beta}(\Sp^{n})} \le C(n,q)  \).

        Assume on the contrary there exist a sequence of solutions \( \sigma_{i} \) such that \( \lambda_{2}(-\Delta_{\Omega_{i}}) \to n \) with corresponding eigenfunctions \( f_{i} \) normalized so that \( ||f_{i}||_{C^{0}} = 1 \). By \( C^{2, \beta} \) estimates we can take a subsequence which is convergent in \( C^{2,\beta'}(\Sp^{n}) \). Relabel the sequence and assume that \( \sigma_{i} \to \sigma_{\infty} \). Now, \( f_{i} \) is orthogonal to \( E_{0}(\Omega_{i}) = \operatorname{span}(1) \) and \( E_{1}(\Omega_{i}) = \operatorname{span}\{z^{l} / \sigma_{i}\} \) in \( L^{2}(\Sp^{n}, dV_{\Omega_{i}}) \) and solves the equation 
        \[ -\Delta_{\Omega_{i}}f_{i} = \lambda_{2}(-\Delta_{\Omega_{i}})f_{i}, \]
        which is a linear second-order elliptic equation with bounded ellipticity (thanks to \( C^{2,\beta} \) estimates). Schauder estimates imply uniform \( C^{2, \beta} \) estimates on \( f_{i} \), so we can take a convergent \( C^{2} \) subsequence with \( f_{i} \to f_{\infty} \). The orthogonality also passes to the limit, so we have 
        \[ f_{\infty} \in (E_{0} \oplus E_{1})^{\perp} \text{ in } L^{2}(\Sp^{n}, dV_{\infty}) \]
        but \( \lambda_{2}(-\Delta_{\Omega_{i}}) \to n \) implies \( -\Delta_{\Omega_{\infty}}f_{\infty} = n f_{\infty} \)
        which forces \( f_{\infty} \in E_{1} \) by the rigidity in \cref{Poincare}. Thus, we arrive at a contradiction.
        
        
    \end{proof}
    \begin{remark}
        The logarithmic Brunn-Minkowski conjecture is equivalent to spectral gap \( \lambda_{2}(-\Delta_{\Omega}) \ge n + 1 \) on even test functions for any even bounded convex body \( \Omega \). An AI proof of the conjecture has been released in \cite{OAI:The-logarithmic-Brunn-Minkowski-conjecture-September-23-2026}. The above proof works only for the solutions of the \( L_{p} \) Minkowski problem but has the advantage of not assuming any symmetry on the convex body.  
    \end{remark}

    \section{Uniqueness near the affine exponent}

    Ivaki-Milman \cite{ivaki2023uniqueness} observed that position vectors are special test functions with closed formulas. They applied this to the spectral Brunn-Minkowski inequality along with few applications of divergence formula with respect to cone-volume measure. It allowed them to give a new proof of uniqueness for \( p \in (-n-1,-1] \). Du \cite{du2025uniqueness} further utilized this idea to prove uniqueness in the supercritical range for strongly symmetric solutions. Du's idea was to subtract translation as well as scaling components of the position vectors and use a spectral gap which was established from the fact that strongly symmetric solutions satisfied \( C^{0} \) estimates. We redo Du's methods with uniform spectral estimate from the previous section.
    
    As before, let \( E_{0} = \operatorname{span}(1) \) and \( E_{1} = \operatorname{span}\{\left< z,B_{l} \right>/ \sigma\} \) be the first two eigenspaces of \( \Delta_{\Omega} \). Unless otherwise stated, all integrals in this section are over \( \Sp^{n} \).
    \begin{lemma}
        Let \( \Omega \) be a solution of the isotropic \( L_{p} \)-Minkowski problem. Let \( X :\Sp^{n} \to \R^{n+1} \) be the inverse Gauss map. Let \( \{B_{l}\}_{l=1}^{n+1} \) be an orthonormal basis of \( \R^{n+1} \). Then, the scalar component of \( X_{l} = \left< X, B_{l} \right> \) in \( L^{2}(\Sp^{n},dV_{\Omega}) \) is given by 
        \[ \pi_{E_{0}}(X_{l}) = \frac{\left< X_{l},1 \right>_{L^{2}(\Sp^{n},dV_{\Omega})}}{\left< 1,1 \right>_{L^{2}(\Sp^{n},dV_{\Omega})}} = \frac{\int X_{l}dV_{\Omega}}{\int dV_{\Omega}} = \frac{n+p+1}{n} \frac{\int \left< \nabla \sigma, B_{l} \right>dV_{\Omega}}{\int dV_{\Omega}}. \]
        Furthermore, the projection satisfies the following norm inequality \begin{equation}
            \sum_{l=1}^{n+1} \left< \pi_{E_{0}}X_{l}, \pi_{E_{0}}X_{l} \right>_{L^{2}(\Sp^{n},dV_{\Omega})} \le \left( \frac{n+p+1}{n} \right)^{2} \int |\nabla \sigma|^{2}dV_{\Omega}.\label{pi0}
        \end{equation}
    \end{lemma}
    \begin{proof}
        At a solution, the cone-volume measure is given by \( dV_{\Omega} = \frac{\sigma \det A[\sigma]d \mu}{n+1} = \frac{\sigma^{p}d \mu}{n+1} \). Also observe that \( \operatorname{div}_{\Sp^{n}}(B_{l}^{T}) = \operatorname{div}_{\Sp^{n}}(B_{l} - z^{l}z) = - nz^{l} \). Using these two facts, we calculate
        \begin{align*}
            \pi_{E_{0}}(X_{l}) & = \frac{\int \left< \sigma z + \nabla \sigma, B_{l} \right>dV_{\Omega}}{\int dV_{\Omega}} \\
            & = \frac{\frac{1}{n+1}\int \sigma^{p+1}\left< z , B_{l} \right> d\mu + \int \left< \nabla \sigma, B_{l} \right>dV_{\Omega}}{\int dV_{\Omega}} \\
            & = \frac{-\frac{1}{n(n+1)}\int \sigma^{p+1} \operatorname{div}(B_{l}^{T})d \mu+ \int \left< \nabla \sigma, B_{l} \right>dV_{\Omega}}{\int dV_{\Omega}} \\
            & = \frac{ \frac{p+1}{n}\int \left< \nabla \sigma, B_{l} \right>dV_{\Omega}+ \int \left< \nabla \sigma, B_{l} \right>dV_{\Omega}}{\int dV_{\Omega}} \\
            & = \frac{n+p+1}{n} \frac{\int \left< \nabla \sigma, B_{l} \right>dV_{\Omega}}{\int dV_{\Omega}}.
        \end{align*}
        Thus, by Cauchy-Schwarz one gets
        \begin{align*}
            \sum_{l=1}^{n+1} \left< \pi_{E_{0}}X_{l}, \pi_{E_{0}}X_{l} \right>_{L^{2}(\Sp^{n},dV_{\Omega})}  & = \left( \frac{n+p+1}{n} \right)^{2} \sum_{l=1}^{n+1} \frac{\left( \int \left< \nabla \sigma,B_{l} \right>dV_{\Omega} \right)^{2}}{\int dV_{\Omega}} \\
            & \le \left( \frac{n+p+1}{n} \right)^{2} \sum_{l=1}^{n+1} \int \left< \nabla \sigma, B_{l} \right>^{2}dV_{\Omega} = \left( \frac{n+p+1}{n} \right)^{2} \int |\nabla \sigma|^{2}dV_{\Omega}.
        \end{align*}
    \end{proof}
    \begin{lemma}
        Under the same hypothesis as in the previous lemma, the translation component of \( X \) satisfies the following norm inequality \begin{equation}
            \sum_{l=1}^{n+1} \left< \pi_{E_{1}}X_{l}, \pi_{E_{1}}X_{l} \right>_{L^{2}(\Sp^{n},dV_{\Omega})} \le  \int \sigma^{2} dV_{\Omega}. \label{pi1}
        \end{equation}
    \end{lemma}
    \begin{proof}
        We can pick an orthonormal basis in \( \R^{n+1} \) such that \( \frac{z^{i}}{\sigma} \)'s are orthogonal in \( L^{2}(\Sp^{n}, dV_{\Omega}) \). Then the projection onto \( E_{1} = \operatorname{span} \{ z^{i} / \sigma\} \) is given by 
        \[ \pi_{E_{1}}(X_{l}) = \sum_{i=1}^{n+1} \frac{ \left< X_{l},z^{i} /\sigma \right>_{L^{2}(\Sp^{n}, dV_{\Omega})}}{\left< z^{i} / \sigma, z^{i} / \sigma \right>_{L^{2}(\Sp^{n},dV_{\Omega})}} \frac{z^{i}}{\sigma} = \sum_{i=1}^{n+1} \frac{\int X_{l}z^{i} / \sigma dV_{\Omega}}{\int (z^{i})^{2}/\sigma^{2}d V_{\Omega}} \frac{z^{i}}{\sigma}.   \]
        For the numerator, observe that \begin{align*}
            \int_{\Sp^{n}} X_{l}z^{i} / \sigma dV_{\Omega} & = \frac{1}{n+1}\int_{\partial \Omega}\left< X_{l}B_{i},\nu \right>d \mathcal{H}^{n} \\
            & = \frac{1}{n+1}\int_{\Omega} \operatorname{div}_{\R^{n+1}}(X_{l}B_{i})dx = \frac{\delta_{il}}{n+1} \int_{\Omega}dx = \frac{\delta_{il}}{n+1} \int_{\Sp^{n}}dV_{\Omega}.
        \end{align*}
        It follows that
        \begin{align*}
            \sum_{l=1}^{n+1} \left< \pi_{E_{1}}X_{l}, \pi_{E_{1}}X_{l} \right>_{L^{2}(\Sp^{n}, dV_{\Omega})} & = \sum_{l=1}^{n+1} \sum_{i=1}^{n+1} \frac{(\int X_{l}z^{i} / \sigma dV_{\Omega})^{2}}{\int (z^{i})^{2}/\sigma^{2}d V_{\Omega}} \\
            & = \left( \frac{1}{n+1}\int dV_{\Omega} \right)^{2} \sum_{i=1}^{n+1} \frac{1}{\int (z^{i})^2 / \sigma^{2} dV_{\Omega}} \\
            & \le \left(\frac{1}{n+1}\int dV_{\Omega}\right)^{2} \sum_{i=1}^{n+1} \frac{(\int(z^{i})^{2}\sigma^{2}dV_{\Omega})}{(\int (z^{i})^{2}dV_{\Omega})^{2}} \\
            & = \sum_{i=1}^{n+1}\int (z^{i})^{2}\sigma^{2}dV_{\Omega} = \int \sigma^{2}dV_{\Omega} 
        \end{align*}
        where we used Cauchy-Schwarz inequality and isotropy from \cref{isotropy}, in the last line.
    \end{proof}
    We now prove the main theorem.
    
    \begin{proof}[Proof of \cref{maintheorem}]
    Let \( \Omega \) be a solution of \cref{LpMinkowski}. To use the eigenvalue estimate, we remove the scaling and translation components. Let \( \tilde{X}_{l} = X_{l} - \pi_{E_{0}}(X_{l}) - \pi_{E_{1}}(X_{l}) \) and note that,
        \begin{align}
             \sum_{l=1}^{n+1}\left< X_{l}, \Delta_{\Omega}X_{l} \right> & =  -\sum_{l=1}^{n+1}\int |\nabla X_{l}|_{g_{\Omega}}^{2}dV_{\Omega} \nonumber \\
            & =  -\sum_{l=1}^{n+1} \int \partial_{i}X_{l}\partial_{j}X_{l}(A^{-1})^{ij}\sigma dV_{\Omega} \nonumber\\
            & =  - \int \sigma (\Delta \sigma + n \sigma)dV_{\Omega} \nonumber \\
            & = -\frac{1}{n+1}\int \sigma^{p+1}\Delta \sigma d \mu - n \int \sigma^{2}dV_{\Omega} = (p+1)\int |\nabla \sigma|^{2}dV_{\Omega} - n \int \sigma^{2}dV_{\Omega}. \label{XDX}
        \end{align}
    See \cite[Lemma 3.2]{ivaki2023uniqueness} or \cite[Lemma 2.2]{du2025uniqueness} for a detailed derivation of this equality. 
    Also, \begin{equation}
        \sum_{l=1}^{n+1} \left< X_{l},X_{l} \right> = \sum_{l=1}^{n+1} \int X_{l}^{2}dV_{\Omega} = \int (\sigma^{2} + |\nabla \sigma|^{2})dV_{\Omega}. \label{XX}
    \end{equation}
    Now, as \( \tilde{X}_{l} \in (E_{0} \oplus E_{1})^{\perp} \), we have \begin{align*}
        \sum_{l=1}^{n+1} -\left< \tilde{X}_{l}, \Delta_{\Omega} \tilde{X}_{l} \right> &\ge \lambda_{2} \sum_{l=1}^{n+1} \left< \tilde{X}_{l}, \tilde{X}_{l} \right>.
    \end{align*}
    Expanding both sides,
    \begin{align*}
        \sum_{l=1}^{n+1} - \left< X_{l},\Delta_{\Omega} X_{l} \right> - n\left< \pi_{E_{1}}X_{l}, \pi_{E_{1}}X_{l} \right> & \ge \lambda_{2} \sum_{l=1}^{n+1} \left( \left< X_{l},X_{l} \right> - \left< \pi_{E_{1}}X_{l}, \pi_{E_{1}}X_{l} \right> - \left< \pi_{E_{0}}X_{l}, \pi_{E_{0}}X_{l} \right> \right).
    \end{align*}
    Rearranging the inequality and using \cref{XDX} and \cref{XX} gives 
    \begin{align*}
        -\int ((p+1)|\nabla \sigma|^{2} - n \sigma^{2})dV_{\Omega} + (\lambda_{2}-n) \sum_{l=1}^{n+1}\left< \pi_{E_{1}}X_{l}, \pi_{E_{1}}X_{l} \right> +  \lambda_{2} \sum_{l=1}^{n+1}\left< \pi_{E_{0}}X_{l}, \pi_{E_{0}}X_{l} \right> & \ge \lambda_{2}\int(\sigma^{2} + |\nabla \sigma|^{2})dV_{\Omega}.
    \end{align*}
    Utilizing the norm inequalities from the previous lemmas and rearranging further gives
    \begin{align}
        \lambda_{2}\left( \frac{n+p+1}{n} \right)^{2} \int |\nabla \sigma|^{2}dV_{\Omega} & \ge \left( \lambda_{2} + 1+ p \right) \int |\nabla \sigma|^{2}dV_{\Omega}. \label{final}
    \end{align}
    Let \( \delta = -p-n-1 > 0 \), we can solve for \( \lambda_{2} + 1+p > \lambda_{2}\left( \frac{n+p+1}{n} \right)^{2} \), or equivalently
    \[ \lambda_{2} -n- \delta - \lambda_{2}\frac{\delta^{2}}{n^{2}} > 0. \]
    From \cref{spectralestimate}, \( \lambda_{2} \ge n + \epsilon \) so we can find a \( \delta_{n}> 0 \) such that for  \( \delta \in (0, \delta_{n}) \) the linear coefficient dominates the quadratic coefficient in inequality \cref{final} and we get \( |\nabla \sigma| =0 \) for the range of \( \delta \), so uniqueness holds. 
    \end{proof}

\bibliographystyle{acm}
\bibliography{ref}
\end{document}